\documentclass[twoside,leqno,twocolumn]{article}
\usepackage{ltexpprt, amsmath}
\usepackage{amssymb}

\newtheorem{thm}{theorem}[section]
\newtheorem{claim}[thm]{Claim}

\newcommand{\qed}[0]{{\hspace*{\fill}\mbox{$\Box$}}}

\newcommand{\ul}[0]{\underline}

\newcommand{\cA}[0]{{\cal A}}

\newcommand{\cC}[0]{{\cal C}}

\newcommand{\cE}[0]{{\cal E}}

\newcommand{\cI}[0]{{\cal I}}

\newcommand{\cM}[0]{{\cal M}}

\newcommand{\cO}[0]{{\cal O}}

\newcommand{\cW}[0]{{\cal W}}

\newcommand{\Z}{{\mathbb Z}}

\newcommand{\gb}[0]{\beta}

\newcommand{\grg}[0]{\gamma}

\newcommand{\go}[0]{\omega}
\newcommand{\gO}[0]{\Omega}
\newcommand{\gs}[0]{\sigma}

\newcommand{\ant}{{\operatorname{\rm int}\,}}
\newcommand{\comp}{{\operatorname{\rm comp}}}
\makeatletter
\def\eqalign#1{\,\vcenter{\openup\jot\m@th
  \ialign{\strut\hfil$\displaystyle{##}$&$\displaystyle{{}##}$\hfil
      \crcr#1\crcr}}\,}
\makeatother

\begin{document}
\renewcommand{\thefootnote}{\fnsymbol{footnote}}

\title{Torpid Mixing of Local Markov Chains on 3-Colorings\\ of the Discrete Torus}

\author{David Galvin\thanks{Department of Mathematics, University of Pennsylvania,
Philadelphia, PA 19104; dgalvin@math.upenn.edu.}\\Dana
Randall\thanks{College of Computing, Georgia Institute of
Technology, Atlanta, GA 30332; randall@cc.gatech.edu.  Supported in part by
NSF grants CCR-0515105 and DMS-0505505.}}

\date{}

\maketitle

\begin{abstract}
We study local Markov chains for sampling 3-colorings of the discrete torus
$T_{L,d}=\{0,\ldots, L-1\}^d$.
We show that there is a constant $\rho \approx .22$ such that for
all even $L \geq 4$ and $d$ sufficiently large, certain local Markov chains require
exponential time to converge to equilibrium.  More precisely,
if $\cM$ is a Markov chain on the set of proper $3$-colorings
of $T_{L,d}$ that updates the color of at most $\rho  L^d$ vertices
at each step and whose stationary distribution is uniform, then the
convergence to stationarity of $\cM$ is exponential in
$L^{d-1}$.
Our proof is based on a conductance argument that builds on sensitive
new combinatorial enumeration techniques.
\end{abstract}

\section{Introduction}
Sampling and counting colorings of a graph are fundamental problems
in computer science and discrete mathematics.  We consider the
problem of sampling uniformly at random from the set
$\cC_k=\cC_k(G)$ of proper $k$-colorings of a graph $G = (V , E)$. A
proper $k$-coloring $\chi$ is a labeling $\chi : V \rightarrow k$
such that all neighboring vertices have different colors. This
sampling problem is also fundamental in statistical physics and
corresponds to generating configurations from  the Gibbs
distribution of the zero-temperature antiferromagnetic Potts model
\cite{welsh}. From the physics perspective, the underlying graph is
typically taken to be the cubic lattice $\Z^d$ and sampling and
counting reveal underlying thermodynamic properties of the
corresponding physical system.

Much focus has gone towards solving the sampling problem using
rapidly mixing Markov chains. The idea is to design a Markov chain
whose stationary distribution is uniform over the set of proper
colorings. Then, starting at an arbitrary coloring and simulating a
random walk according to this chain for a sufficient number of
steps, we get a sample from close to the desired distribution. The
number of steps required of this walk is referred to as the {\it
mixing time}   (see, e.g., \cite{sinc}). The chain is called {\it
rapidly mixing} if the mixing time is polynomial in $n = |V |$ (so
it converges quickly to stationarity); it is {\it torpidly mixing}
if its mixing  time is super-polynomial in $n$ (so it converges
slowly). There has been a long history of studying mixing times of
various chains in the context of colorings (see, e.g., \cite{ammb,
FriezeVigoda, gmp, hv, jer, lrs}).

A particular focus of this study has been on {\em Glauber dynamics}.
For proper $k$-colorings this is any single-site update Markov chain
that connects two colorings only if they differ on at most a single
vertex. The {\it Metropolis} chain $\cM_k$ on state space $\cC_k$
has transition probabilities $P_k(\chi_1,\chi_2)$, $\chi_1,\chi_2
\in \cC_k,$ given by
$$
P_k(\chi_1,\chi_2) = \left\{
            \begin{array}{ll}
               \! \!  0, & \! \! \! \! \!  \mbox{ if $|\{v \in V:\chi_1(v)\neq \chi_2(v)\}| > 1$}; \\
               & \\
               \! \! \frac{1}{k|V|},  &  \! \! \! \! \! \mbox{ if $|\{v \in V:\chi_1(v)\neq \chi_2(v)\}| = 1$}; \\
               & \\
               \! \! 1 -  &
              \! \! \! \! \! \! \! \! \!
               \sum_{\chi_1 \neq \chi_2' \in \cC_k}  P_k(\chi_1,\chi_2'), \\
& \! \! \! \! \! \mbox{ if $\chi_1=\chi_2$.}
            \end{array}
         \right.
$$
We may think of $\cM_k$ dynamically as follows. From a $k$-coloring
$\chi$, choose a vertex $v$ uniformly from $V$ and a color $j$
uniformly from $\{0,\ldots,k-1\}$.  Then recolor $v$ with color $j$
if this is a proper $k$-coloring; otherwise stay at $\chi$.

When $\cM_k$ is ergodic, its stationary distribution $\pi_k$ is
uniform over proper $k$-colorings. A series of recent papers have
shown that  $\cM_k$ is rapidly mixing provided the number of colors
is sufficiently large compared to the maximum degree (see
\cite{FriezeVigoda} and the references therein). Substantially less
is known when the number of colors is small.  In fact, for $k$ small
it is NP-complete to decide whether a graph admits even one
$k$-coloring.

In this paper we focus on the mixing rate of $\cM_k$ on rectangular
regions of the cubic lattice ${\mathbb Z}^d$. Observe that the
lattice is  bipartite, so it always admits a $k$-coloring for any
$k\geq 2$.   It is also known that Glauber dynamics connects the
state space of $k$-colorings on any such lattice region \cite{lrs}.
In $\Z^2$ much is known about the mixing rate of $\cM_k$. Luby et
al. \cite{lrs} showed that Glauber dynamics for sampling 3-colorings
is rapidly mixing on any finite, simply-connected subregion of
$\Z^2$ when the colors on the boundary of the region are fixed.
Goldberg et al. \cite{gmp} subsequently showed that the chain
remains fast on rectangular regions without this boundary
restriction. Substantially more is known when there are many colors:
Jerrum \cite{jer} showed that Glauber dynamics is rapidly mixing on
any graph satisfying $k \geq 2 \Delta$, where $k$ is the number of
colors and $\Delta$ is the maximum degree, thus showing Glauber
dynamics is fast on $\Z^2$ when $k \geq 8$.  It has since been shown
that it is fast for $k \geq 6$ \cite{ammb, bdg}.  Surprisingly the
efficiency remains unresolved for $k=4$ or~$5$.

In higher dimensions much less is known when $k$ is small.
Physicists have performed extensive numerical experiments \cite{fs,
swk} suggesting that Glauber dynamics on 3-colorings is torpidly
mixing when the dimension of the cubic lattice is large enough. We
prove this conjecture for the first time here by studying the mixing
time of the chain on cubic lattices with periodic boundary
conditions.

\subsection{Results}
Our focus in this paper is sampling 3-colorings of the even discrete
torus $T_{L,d}$. This is the graph on vertex set $\{0, \ldots,
L-1\}^d$ (with $L$ even) with edge set consisting of those pairs of
vertices that differ on exactly one coordinate and differ by $1$
(mod $L$) on that coordinate. For a Markov chain $\cM$ on the
3-colorings of $T_{L,d}$ we denote by $\tau_\cM$ the mixing time of
the chain; this will be formally defined in Section
\ref{sec-partitioning}. Our main theorem is the following.
\begin{thm} \label{thm-glauber.slow.mixing}
There is a constant $d_0>0$ for which the following holds. For
$d\geq d_0$ and $L\geq 4$ even, the Glauber dynamics chain $\cM_3$
on $\cC_3(T_{L,d})$ satisfies
$$
\tau_{\cM_3} \geq \exp \left\{\frac{L^{d-1}}{d^4\log^2 L}\right\}.
$$
\end{thm}

Our techniques actually apply to a more general class of chains. A
Markov chain $\cM$ on state space $\cC_3$ is {\em $\rho$-local} \
if, in each step of the chain, at most $\rho|V|$ vertices have
their color changed; that is, if
$$
P_\cM(\chi_1,\chi_2) \neq 0$$
implies
$$ |\{v \in V : \chi_1(v)\neq
\chi_2(v)\}| \leq \rho |V|.
$$
These types of chains were introduced in \cite{DyerFriezeJerrum},
where the terminology {\em $\rho|V|$-cautious} was employed. We
prove the following, which easily implies
Theorem~\ref{thm-glauber.slow.mixing}.
\begin{thm} \label{thm-rho.local.slow.mixing}
Fix $\rho>0$ satisfying $H(\rho)+\rho < 1$. There is a constant
$d_0=d_0(\rho)>0$ for which the following holds. For $d\geq d_0$ and
$L\geq 4$ even, if $\cM$ is an ergodic $\rho$-local Markov chain on
$\cC_3(T_{L,d})$ with uniform stationary distribution then
$$
\tau_{\cM} \geq \exp \left\{\frac{L^{d-1}}{d^4\log^2 L}\right\}.
$$
\end{thm}
Here $H(x)=-x\log x -(1-x)\log (1-x)$ is the usual binary entropy
function. Note that all $\rho \leq .22$ satisfy $H(\rho)+\rho < 1$.

\subsection{Techniques} We show slow mixing via a
conductance argument by identifying a ``bad cut'' in the state space
requiring exponential time to cross. Intuitively, in sufficiently
high dimension, the set of 3-colorings of the lattice is believed to
naturally partition into 6 classes: each class is identified by a
predominance of one (of 3) colors on one of the two (even or odd)
sublattices.  This characterization was recently rigorously verified
on the infinite lattice, thereby establishing the existence of 6
distinct ``maximal entropy Gibbs states''  \cite{GKRS}.   That work
builds heavily on technical machinery introduced by Galvin and Kahn
\cite{GalvinKahn} showing that independent sets partition similarly
in sufficiently high dimensions in that they lie primarily on the
even or odd sublattices. Specifically, write $\cE$ and $\cO$ for the
sets of even and odd vertices of ${\mathbb Z}^d$ (defined in the
obvious way) and set $\Lambda_L =[-L,L]^d$ and $\partial \Lambda_L
=[-L,L]^d \setminus [-(L-1),L-1]^d$. For $\lambda>0$, choose
${\mathbb I}$ from ${\cal I}(\Lambda_L)$ (the set of independent
sets of the box) with $\Pr({\mathbb I}=I) \propto \lambda^{|I|}$.
Galvin and Kahn showed that for $\lambda > Cd^{-1/4}\log^{3/4}d$
(for a large constant $C$) and fixed $v \in \Lambda_L \cap \cE$
$$
\eqalign{ \lim_{L\rightarrow\infty}{\mathbb P}(v\in{\mathbb I} ~  |
& ~  {\mathbb I}  \supseteq
\partial \Lambda_L  \cap {\cal E}) \cr
& >  \lim_{L\rightarrow\infty}{\mathbb P}\left(v\in{\mathbb I}~|~{\mathbb
I}\supseteq
\partial \Lambda_L\cap {\cal O}\right). }
$$
In other words, the influence of the boundary on the center of a
large box persists as the boundary recedes.

Notice that neither of the results of \cite{GalvinKahn} or
\cite{GKRS} establishing the presence of multiple Gibbs states
directly implies anything about the behavior of Markov chains on
finite lattice regions. However, they do suggest that in the finite
setting, typical configurations fall into the distinct classes
described in stationarity and that it will be unlikely to move
between these classes; the remaining configurations are expected to
have negligible weight for large lattice regions, even when they are
finite.

Galvin \cite{Galvin} extended the results of \cite{GalvinKahn},
showing that in sufficiently high dimension, Glauber dynamics on
independent sets mixes slowly in rectangular regions of ${\mathbb
Z}^d$ with periodic boundary conditions. Similar results were known
previously about independent sets; however, one significant new
contribution of \cite{Galvin} was showing that as $d$ increases, the
critical $\lambda$ above which Glauber dynamics mixes slowly tends
to $0$. In particular, there is some dimension $d_0$ such that for
all $d \geq d_0$, Glauber dynamics will be slow on $\Z^d$ when
$\lambda=1$. This turns out to be the crucial new ingredient
allowing us to rigorously verify slow mixing for sampling
3-colorings in high dimensions, as there turns out to be a close
connection between the independent set model at $\lambda=1$ and the
$3$-coloring model. Note that unlike most statistical physics
models, the 3-coloring problem does not have a parameter $\lambda$
that can be tweaked to establish desired bounds; this makes the
proofs here significantly more delicate than the usual slow mixing
arguments. Section \ref{subsec-overview} provides a more detailed
discussion of the elements of the proof and of some of the
difficulties inherent to the sampling problem under discussion.

\section{Partitioning the state space} \label{sec-partitioning}
We begin by formalizing some definitions. Given a Markov chain
${\cal M}$ on state space $\Omega$ with uniform stationary
distribution denoted by $\pi$, let $P^t(X,
\cdot)$ be the distribution of the chain at time $t$ given that it
started in state $X$.
The {\em mixing time} $\tau_{{\cal M}}$ of $\cM$ is defined to be
$$
\tau_{\cM}=\min \left\{t_0:  ||P^t, \pi||_{\mbox{tv}} \leq
\frac{1}{e} ~~~ \forall t>t_0\right\}
$$
where
$$ ||P^t, \pi||_{\mbox{tv}} = \max_{X \in \Omega} \frac{1}{2} \sum_{Y \in \Omega}
|P^t(X, Y) - \pi(Y)|,$$
is the {\it total variation distance}.


We prove Theorem \ref{thm-rho.local.slow.mixing} via a well-known
{\it conductance argument} \cite{JerrumSinclair, LawlerSokal, Thomas}, using
a form of the argument derived in \cite{DyerFriezeJerrum}.
As above, let $\cM$ be an ergodic Markov chain on state space $\gO$ with
transition probabilities $P$ and stationary distribution $\pi$. Let
$A \subseteq \gO$ and $M \subseteq \gO \setminus A$ satisfy $\pi(A)
\leq 1/2$ and $\go_1 \in A, \go_2 \in \gO \setminus (A \cup M)
\Rightarrow P(\go_1, \go_2) =0$. Then from \cite{DyerFriezeJerrum}
we have
\begin{equation} \label{conductance_bound}
\tau_\cM \geq \frac{\pi(A)}{8\pi(M)}.
\end{equation}

Let us return to the setup of Theorem
\ref{thm-rho.local.slow.mixing}. For even $L$, $T_{L,d}$ is
bipartite with partition classes $\cE$ (consisting of those vertices
the sum of whose coordinates is even) and $\cO$. To show torpid
mixing, it is sufficient to identify a single bad cut.  We
concentrate on the vertices in each 3-coloring that are colored with
the first color, 0.  The objective of
Theorem~\ref{thm-rho.local.slow.mixing} will be to verify that most
3-colorings have an imbalance whereby the vertices colored 0 lie
predominantly on $\cE$ or $\cO$, and those that are roughly balanced
on the two sublattices are highly unlikely in stationarity.  This is
sufficient to show that the conductance is small.

Accordingly let us define the set of ``balanced'' $3$-colorings by
$$
\cC_3^{b,\rho} = \{\chi \in \cC_3:\left||\chi^{-1}(0)\cap
\cE|\!-\!|\chi^{-1}(0)\cap \cO|\right| \leq \rho L^d/2\}
$$
and likewise let
$$
\cC_3^{\cE,\rho} = \{\chi \in \cC_3:|\chi^{-1}(0)\cap
\cE|>|\chi^{-1}(0)\cap \cO| + \rho L^d/2\}.
$$
By symmetry,
$\pi_3(\cC_3^{\cE,\rho}) \leq 1/2$. Notice that since $\cM$
updates at most $\rho L^d$ vertices in each step, we have that
if $\chi_1 \in \cC_3^{\cE,\rho}$ and $\chi_2 \in \cC_3 \setminus
(\cC_3^{\cE,\rho} \cup \cC_3^{b,\rho})$ then $P_\cM(\chi_1,\chi_2) =
0$.
Therefore, by (\ref{conductance_bound}),
$$
\tau_\cM \geq \frac{\pi_3(\cC_3^{\cE,\rho})}{8\pi_3(\cC_3^{b,\rho})}
\geq \frac{1-\pi_3(\cC_3^{\cE,\rho})}{16\pi_3(\cC_3^{b,\rho})},
$$
and so Theorem \ref{thm-rho.local.slow.mixing} follows from the
following critical theorem.
\begin{thm} \label{thm-main}
Fix $\rho>0$ satisfying $H(\rho)+\rho < 1$. There is a constant
$d_0=d_0(\rho)>0$ for which the following holds. For $d \geq d_0$
and $L \geq 4$ even,
$$
\pi_3(\cC_3^{b,\rho}) \leq \exp\left\{\frac{-2L^{d-1}}{d^4\log^2
L}\right\}.
$$
\end{thm}

\section{Proof of Theorem \ref{thm-main}} \label{sec-proof.of.mn.thm}

\subsection{Setup and overview} \label{subsec-overview}

For a generic $\chi \in \cC_3^{b,\rho}$ there are regions of
$T_{L,d}$ consisting predominantly of even vertices colored $0$
together with their neighbors, and regions consisting of odd
vertices colored $0$ together with their neighbors. These regions
are separated by
two-layer ``$0$-free'' moats or {\em cutsets}. In Section
\ref{subsec-cutsets} we describe a procedure that selects a
particular collection of these cutsets. Our main technical result,
Lemma \ref{lem-volume.bounds.from.gk}, asserts that for each
specification of cutset sizes $c_1, \ldots, c_\ell$ and vertices
$v_1, \ldots, v_\ell$, the probability that a coloring has among its
associated cutsets a collection $\grg_1, \ldots, \grg_\ell$ with
$|\grg_i|=c_i$ and with $v_i$ surrounded by $\grg_i$ is
exponentially small in the sum of the $c_i$'s. This lemma is
presented in Section \ref{subsec-results.from.GalvinKahn} and
Theorem \ref{thm-main} is derived from it in Section
\ref{subsec-proof_of_thm}.

The main thrust of \cite{GKRS} is the proof of a result that is
essentially (but not quite) the case $\ell=1$ of Lemma
\ref{lem-volume.bounds.from.gk}. One difficulty we have to overcome
in moving from a Gibbs measure argument to a torpid mixing argument
is that of going from bounding the probability of a configuration
having a single cutset to bounding the probability of it having an
ensemble of cutsets. Another difficulty is that the cutsets we
consider in these ensembles can be topologically more complex than
the connected cutsets that are considered in \cite{GKRS}. In part,
both of these difficulties are dealt with by the machinery developed
in \cite{Galvin}.

We use a ``Peierl's argument'' to prove Lemma
\ref{lem-volume.bounds.from.gk}. By carefully modifying each $\chi
\in \cC_3^{b,\rho}$ inside its cutsets,
we can exploit the fact that the cutsets are $0$-free to map $\chi$
to a set $\varphi(\chi)$ of many different $\chi' \in \cC_3$.
If the $\varphi(\chi)$'s were disjoint for distinct $\chi$'s, we
would essentially be done, having shown that there are many more
$3$-colorings in total than $3$-colorings in $\cC_3^{b,\rho}$.
To control the possible overlap,
we define a flow $\nu:\cC_3^{b,\rho} \times\cC_3\rightarrow
[0,\infty)$ supported on pairs $(\chi,\chi')$ with $\chi' \in
\varphi(\chi)$ in such a way that the flow out of each $\chi \in
\cC_3^{b,\rho}$ is $1$. Any uniform bound we can obtain on the flow
into elements of $\cC_3$
is then easily seen to be a bound on $\pi_3(\cC_3^{b,\rho})$. We
define the flow via a notion of approximation modified from
\cite{GalvinKahn}. To each cutset $\grg$ we associate a set
$A(\grg)$ that approximates the interior of $\grg$ in a precise
sense, in such a way that as we run over all possible $\grg$, the
total number of approximate sets used is small. Then for each $\chi'
\in \cC_3$ and each collection of approximations $A_1, \ldots,
A_\ell$, we consider the set of those $\chi \in \cC_3^{b,\rho}$ with
$\chi' \in \varphi(\chi)$ and with $A_i$ the approximation to
$\grg_i$. We define the flow so that if this set is large, then
$\nu(\chi,\chi')$ is small for each $\chi$ in the set. In this way
we control the flow into $\chi'$ corresponding to each collection of
approximations $A_1, \ldots, A_\ell$; since the total number of
approximations is small, we control the total flow into $\chi'$. In
the language of statistical physics,
this approximation scheme
is a {\em course-graining} argument. The details appear in
Section~\ref{sec-proofs}.

The main results of \cite{Galvin} and \cite{GalvinKahn} are proved
along similar lines to those described above. One of the
difficulties we encounter in moving from these arguments on
independent sets to arguments on colorings is that of finding an
analogous way of modifying a coloring inside a cutset in order to
exploit the fact that it is $0$-free. The beginning of Section
\ref{sec-proofs} (in particular Claims \ref{claim-shift.works} and
\ref{claim-unique.reconstruction}) describes an appropriate
modification that has all the properties we desire.

\subsection{Cutsets} \label{subsec-cutsets}

We describe a way of associating with each $\chi \in \cC_3^{b,\rho}$
a collection of minimal edge cutsets, following the approaches of
\cite{BorgsChayesFriezeKimTetaliVigodaVu} and \cite{Galvin}. First
we need a little notation.

Write $V$ for the vertex set of $T_{L,d}$ and $E$ for its edge set.
For $X \subseteq V$, write $\nabla(X)$ for the set of edges in $E$
that have one end in $X$ and one end outside $X$; $\overline{X}$ for
$V \setminus X$; $\partial_{int}X$ for the set of vertices in $X$
that are adjacent to something outside $X$; $\partial_{ext}X$ for
the set of vertices outside $X$ that are adjacent to something in
$X$; $X^+$ for $X \cup
\partial_{ext} X$; $X^\cE$ for $X \cap \cE$ and $X^\cO$ for $X \cap
\cO$. Further, for $x \in V$ set $\partial x=\partial_{ext}\{x\}$.
We abuse notation slightly, identifying
sets of vertices of $V$ and the subgraphs they induce.


For each $\chi \in \cC_3^{b,\rho}$ set $I=I(\chi)=\chi^{-1}(0)$.
Note that $I(\chi)$ is an independent set (a set of vertices no two
of which are adjacent). For each component $R$ of $(I^\cE)^+$ or
$(I^\cO)^+$ and each component $C$ of $\overline{R}$, set $\gamma =
\gamma_{RC}(I)=\nabla(C)$ and $W=W_{RC}(I)=\overline{C}$. Evidently
$C$ is connected, and $W$ consists of $R$, which is connected,
together with a number of other components of $\overline{R}$, each
of which is connected and joined to $R$, so $W$ is connected also.
It follows that $\grg$ is a minimal edge-cutset in $T_{L,d}$. Say
that $\grg$ is {\em even} if $R$ is a component of $(I^\cE)^+$
and {\em odd} otherwise. Define $\ant \gamma$, the {\em interior} of
$\gamma$, to be the smaller of $C,W$ (if $|W|=|C|$, take $\ant
\gamma = W$).

The cutsets $\grg$ associated to $\chi$ depend
only on the independent set $I(\chi)$, and coincide exactly with the
cutsets associated to an independent set in \cite{Galvin}. We may
therefore apply the machinery developed in \cite{Galvin} for
independent set cutsets in the present setting. In particular, from
\cite[Lemmas 3.1 and 3.2]{Galvin} we know that for each $\chi \in
\cC_3$ there is a collection of associated cutsets $\Gamma(I)$
such that either
\begin{equation} \label{cutset.conditions}
\begin{array}{c}
\mbox{for all $\grg,\grg' \in \Gamma(I)$,} \\
\mbox{$\grg,\grg'$ are even with $\ant \grg \cap \ant \grg' = \emptyset$,} \\
\mbox{and  $I^\cE \subseteq \cup_{\grg} \ant \grg$},
\end{array}
\end{equation}
or we have the analogue of (\ref{cutset.conditions}) with even replaced by
odd. Set ${\cC}_3^{even} = \{\chi \in \cC_3:\chi~\mbox{satisfies
(\ref{cutset.conditions})}\}$. From here on whenever $\chi \in
{\cC}_3^{even}$ is given we assume that $I$ is its associated
independent set and that $\Gamma(I)$ is a particular collection of
cutsets associated with $\chi$ and satisfying
(\ref{cutset.conditions}).
Numerous properties of $\gamma \in \Gamma(I)$ are established in
\cite[Lemmas 3.3 and 3.4]{Galvin}. We list some here that will be of
use in the sequel. That the cutsets are indeed $0$-free regions is
established by (\ref{contour.prop.1}).
\begin{equation} \label{contour.prop.0}
\partial_{int}W \subseteq \cO~~\mbox{and}~~\partial_{ext}W \subseteq
\cE;
\end{equation}
\begin{equation} \label{contour.prop.1}
\partial_{int}W \cap I =
\emptyset~~\mbox{and}~~\partial_{ext}W \cap I = \emptyset;
\end{equation}
\begin{equation} \label{contour.prop.2}
W^\cO = \partial_{ext} W^\cE~~\mbox{and}~~W^\cE=\left\{y \in
\cE:\partial y \subseteq W^\cO\right\};
\end{equation}
\begin{equation} \label{contour.prop.3}
\mbox{for large enough $d$, ~~$|\grg| \geq \max\{|W|^{1-1/d},
d^{1.9}\}$}.
\end{equation}

\subsection{The main lemma} \label{subsec-results.from.GalvinKahn}

For $c \in {\mathbb N}$ and $v \in V$ set
$$
\cW(c,v) = \left\{\grg : \begin{array}{l}
\mbox{$\grg \in \Gamma(I)$ for some $\chi \in \cC_3^{even}$} \\
\mbox{with $|\grg|=c$, $v \in W^\cE$}
\end{array}
\right\}
$$
and set $\cW = \cup_{c,v} \cW(c,v)$. A {\em profile} of a collection
$\{\grg_1, \ldots, \grg_\ell\} \subseteq \cW$ is a vector
$\ul{p}=(c_1 ,v_1, \ldots, c_\ell, v_\ell)$ with $\grg_i \in
\cW(c_i,v_i)$ for all $i$. Given a profile vector $\ul{p}$ set
$$
\cC_3(\ul{p}) = \left\{\chi \in \cC_3^{even} : \begin{array}{l}
\mbox{$\Gamma(I)$ contains a subset}\\
\mbox{with profile $\ul{p}$} \end{array}\right\}.
$$
Our main lemma (c.f. \cite[Lemma 3.5]{Galvin})
is the following.
\begin{lemma} \label{lem-volume.bounds.from.gk}
There are constants $c, d_0>0$ such that the following holds. For
all even $L\geq 4$, $d \geq d_0$ and profile vector $\ul{p}$,
\begin{equation} \label{inq-bound.on.inner}
\pi_3(\cC_3(\ul{p})) \leq \exp\left\{-\frac{c\sum_{i=1}^\ell c_i}{d}
\right\}.
\end{equation}
\end{lemma}

We will derive Theorem \ref{thm-main} from Lemma
\ref{lem-volume.bounds.from.gk} in Section \ref{subsec-proof_of_thm}
before proving the lemma in Section \ref{sec-proofs}. From here on
we assume that the conditions of Theorem \ref{thm-main} and Lemma
\ref{lem-volume.bounds.from.gk} are satisfied (with $d_0$
sufficiently large to support our assertions).

\subsection{Proof of Theorem \ref{thm-main} assuming the main lemma} \label{subsec-proof_of_thm}

We begin with an easy count that dispenses with colorings where
$|I(\chi)|$ is small. Set
$$
\cC_3^{small} = \left\{\chi \in
\cC_3^{b,\rho}:\min\{|I^\cE|,|I^\cO|\} \leq L^d/4d^{1/2}\right\}.
$$
\begin{lemma} \label{lem-bounding_small} $\pi_3(\cC_3^{small}) \leq
\exp\left\{-\Omega(L^d)\right\}$.
\end{lemma}

\noindent {\em Proof: }
For any $A \subseteq \cE$ and $B \subseteq \cO$, let
$\comp(A,B)$ be the number of components in
$V \setminus (A \cup B \cup \partial^\star A \cup \partial^\star B)$, where
for $T \subseteq \cE$ (or $\cO$),
$$
\partial^\star T = \{x \in \partial_{ext} T : \partial x \subseteq T\}~(=\{x \in V :
\partial x \subseteq T\}).
$$
We begin by noting that by $\cE$-$\cO$
symmetry
\begin{equation} \label{inq-small.count}
|\cC_3^{small}| \leq 2 \sum 2^{|\partial^\star A|+|\partial^\star
B|+\comp(A,B)},
\end{equation}
where the sum is over all pairs $A \subseteq \cE$, $B \subseteq \cO$
with no edges between $A$ and $B$ and satisfying $|A|\leq
L^d/4d^{1/2}$ and $|B| \leq (\rho + 1/2d^{1/2})L^d/2$.
Indeed, once
we have specified that the set of vertices colored $0$ is $A \cup
B$, we have a free choice between $1$ and $2$ for the color at $x
\in  \partial^\star A \cup
\partial^\star B$, and we also have a free choice between the two
possible colorings of each component of $V \setminus (A \cup B \cup
\partial^\star A \cup \partial^\star B)$.

A key observation is the following.  For $A$ and $B$ contributing to the sum in
(\ref{inq-small.count}),
\begin{equation} \label{inq-number.of.remaining.comps}
\comp(A,B) \leq L^d/2d.
\end{equation}
To see this, let $C$ be a component of $V \setminus (A \cup B)$. If
$C=\{v\}$ consists of a single vertex, then (depending on the parity
of $v$) we have either $\partial v \subseteq A$ or $\partial v
\subseteq B$ and so $v \in \partial^\star A \cup \partial^\star B$.
Otherwise, let $vw$ be an edge of $C$ with $v \in \cE$ (and so $w
\in \cO$). If $v$ has $k$ edges to $B$ and $u$ has $\ell$ to $A$,
then (since there are no edges from $A$ to $B$)
we have $(k-1)+(\ell-1)\leq 2d-2$ or $k+\ell \leq 2d$. (Here we are
using that in $T_{L,d}$, if $uv \in E$ then
there is a matching between all but one of the
neighbors of $u$ and $v$.)
Since $v$ has $2d-1-k$ edges to $\cO \setminus (B \cup \{w\})$ and
$w$ has $2d-1-\ell$ edges to $\cE \setminus (A \cup \{v\})$ we have
that $|C| = 4d -(k+\ell) \geq 2d$. From this,
(\ref{inq-number.of.remaining.comps}) follows.

Inserting (\ref{inq-number.of.remaining.comps}) into
(\ref{inq-small.count}) and bounding $|\partial^\star A|$ and
$|\partial^\star B|$ by the maximum values of $|A|$ and $|B|$ (valid
since $T \subseteq \cE$ (or $\cO$) satisfies $|T| \leq
|\partial_{ext}T|$, so $|\partial^\star T| \leq |T|$) and with the
remaining inequalities justified below, we have
\begin{eqnarray}
|\cC_3^{small}| & \leq & 2^{\frac{L^d}{2}\left(\rho +
\frac{1}{d^{1/2}} + \frac{1}{d}\right)} \cdot \sum_{i
\leq L^d/4d^{1/2}} {L^d/2 \choose i} \nonumber \\
& & \hspace{.9in} \cdot
\sum_{j \leq (\rho + 1/2d^{1/2})L^d/2} {L^d/2 \choose j} \nonumber \\
& \leq & 2^{\frac{L^d}{2}\left(\rho + \frac{1}{d^{1/2}} + \frac{1}{d} + H\left(\frac{1}{2d^{1/2}}\right) + H\left(\rho + \frac{1}{2d^{1/2}}\right)\right)} \label{int4} \\
& \leq & 2^{\frac{L^d}{2}\left(1-\Omega(1)\right)} \label{int5}
\end{eqnarray}
for sufficiently large $d=d(\rho)$. In (\ref{int4}) we use the
Chernoff bound
$\sum_{i=0}^{[\gb M]}{M\choose i} \leq
2^{H(\gb)M}$ for $\gb \leq \frac{1}{2}$; in (\ref{int5}) we use
$H(\rho)+\rho < 1$. Using $2^{L^d/2} \leq |\cC_3|$,
the lemma follows. \qed

We now consider
$$
\cC_3^{large, even}:=(\cC_3^{b,\rho} \setminus \cC_3^{small}) \cap
\cC_3^{even}.
$$
By Lemma \ref{lem-bounding_small} and $\cE$-$\cO$ symmetry, Theorem
\ref{thm-main} reduces to bounding (say)
\begin{equation} \label{inq-remaining}
\pi_3(\cC_3^{large, even}) \leq
\exp\left\{-\frac{3L^{d-1}}{d^4\log^2 L}\right\}.
\end{equation}

Let $\cC_3^{large, even, nt}$ be the set of
$\chi \in \cC_3^{large, even}$ such that there is a
$\grg \in \Gamma(I)$ with $|\grg|\geq L^{d-1}$ (we think of such
cutsets as being topologically non-trivial (``nt''); see
\cite{Galvin} for an explanation of this) and also let
$\cC_3^{large, even, triv} = \cC_3^{large, even} \setminus
\cC_3^{large, even, nt}$. We assert that
\begin{equation}\label{large2}
\pi_3(\cC_3^{large, even, nt}) \leq
\exp\left\{-\Omega\left(\frac{L^{d-1}}{d}\right)\right\}
\end{equation}
and
\begin{equation}\label{large3}
\pi_3(\cC_3^{large, even, triv})  \leq
\exp\left\{-\frac{4L^{d-1}}{d^4 \log^2 L}\right\};
\end{equation}
this gives (\ref{inq-remaining}) and so completes the proof of
Theorem \ref{thm-main}. Both (\ref{large2}) and (\ref{large3}) are
corollaries of Lemma \ref{lem-volume.bounds.from.gk}, and the steps
are identical to those that are used to bound the measures of
``$\cI_{large, even}^{non-trivial}$'' and ``$\cI_{large, even}^{trivial}$''
in \cite[Section 3.3]{Galvin}.

With the sum below running over all vectors $\ul{p}$ of the form
$(c,v)$ with $v \in V$ and $c \geq L^{d-1}$, and with the
inequalities justified below, we have
\begin{eqnarray}\label{large1}
\eqalign{ \pi_3(\cC_3^{large, even, nt}) ~ &  \leq  ~ \sum_{\ul{p}}
\pi_3(\cC_3(\ul{p}))\cr & \leq ~
L^{2d}\exp\left\{-\Omega\left(\frac{L^{d-1}}{d}\right)\right\} \cr &
\leq ~ \exp\left\{-\Omega\left(\frac{L^{d-1}}{d}\right)\right\}, }
\end{eqnarray}
giving (\ref{large2}). We use Lemma \ref{lem-volume.bounds.from.gk}
in (\ref{large1}). The factor of $L^{2d}$ is for the choices of $c$
and $v$.

The verification of (\ref{large3}) involves finding an $i \in
[\Omega(\log d), O(d\log L)]$ and a set $\Gamma_i(I) \subseteq
\Gamma(I)$ of cutsets with the properties that $|\Gamma_i(I)|
\approx L^d/2^i$, $|\grg| \approx 2^i$ for each $\grg \in
\Gamma_i(I)$ and $\sum_{\grg \in \Gamma_i(I)} |\grg| \approx
L^{d-1}$. The measure of $\cC_3^{large,even,triv}$ is then at most
the product of a term that is exponentially small in $L^{d-1}$ (from
Lemma \ref{lem-volume.bounds.from.gk}), a term corresponding to the
choice of a fixed vertex in each of the interiors, and a term
corresponding to the choice of the collection of lengths. The second
term will be negligible because $\Gamma_i(I)$ is small and the third
will be negligible because all $\grg \in \Gamma_i(I)$ have similar
lengths.


More precisely, for $\chi \in \cC_3^{large,even,triv}$ and $\grg \in
\Gamma(I)$ we have $|\gamma| \geq |\ant \grg|^{1-1/d}$ (by
(\ref{contour.prop.3}))
and
so
$$
\sum_{\grg \in \Gamma(I)} |\gamma|^{d/(d-1)} \geq \sum_{\grg \in
\Gamma(I)} |\ant \gamma| \geq |I^\cE| \geq L^d/4d^{1/2}.
$$
The second inequality is from (\ref{cutset.conditions}) and the
third follows since $\chi \not \in \cC_3^{small}$.

Set $\Gamma_i(I) = \{\grg \in \Gamma(I):2^{i-1} \leq |\grg| <
2^i\}$. Note that $\Gamma_i(I)$ is empty for $2^i < d^{1.9}$ (again
by (\ref{contour.prop.3}))
and for $2^{i-1} > L^{d-1}$ so we may assume
that
\begin{equation} \label{inner.property.1}
1.9 \log d \leq i \leq (d-1)\log L+1.
\end{equation}
Since $\sum_{m=1}^\infty 1/m^2 = \pi^2/6$, there is an $i$ such that
\begin{equation} \label{quad.scale}
\sum_{\grg \in \Gamma_i(I)} |\gamma|^{\frac{d}{d-1}} \geq
\Omega\left(\frac{L^d}{d^{1/2}i^2}\right).
\end{equation}
Choose the smallest such $i$ set $\ell=|\Gamma_i(I)|$. We have
$\sum_{\grg \in \Gamma_i(I)} |\grg| \geq \Omega(\ell 2^i)$ (this
follows from the fact that each $\grg \in \Gamma_i(I)$ satisfies
$|\grg| \geq 2^{i-1}$) and
\begin{equation} \label{inner.property.2}
O\left(\frac{dL^d}{2^i}\right) \geq \ell \geq
\Omega\left(\frac{L^d}{2^{\frac{id}{d-1}}i^2d^{1/2}}\right).
\end{equation}
The first inequality follows from that fact that $\sum_\grg |\grg|
\leq dL^d=|E|$; the second follows from (\ref{quad.scale}) and the
fact that each $\grg$ has $|\grg|^{d/(d-1)} \leq 2^{di/(d-1)}$. We
therefore have $\chi \in \cC_3(\ul{p})$ for some $\ul{p}=(c_1, v_1,
\ldots, c_\ell, v_\ell)$ with $\ell$ satisfying
(\ref{inner.property.2}), with
\begin{equation} \label{inner.property.3}
\sum_{j=1}^\ell c_j \geq O(\ell 2^i),
\end{equation}
with
\begin{equation} \label{inner.property.4}
c_j \leq 2^i
\end{equation}
for each $j$ and with $i$ satisfying (\ref{inner.property.1}). With
the sum below running over all profile vectors $\ul{p}$ satisfying
(\ref{inner.property.1}), (\ref{inner.property.2}),
(\ref{inner.property.3}) and (\ref{inner.property.4}) we have
\begin{eqnarray}
\pi_3(\cC_3^{large,even,triv}) & \leq & \sum_{\ul{p}}
\pi_3(\cC_3(\ul{p})). \label{largeeventriv}
\end{eqnarray}
The right-hand side of (\ref{largeeventriv}) is, by Lemma
\ref{lem-volume.bounds.from.gk}, at most
$$
d\log L~\max_{i~\mbox{satisfying (\ref{inner.property.1})}}2^{\ell
i}{L^d \choose \ell}\exp\left\{-\Omega\left(\frac{\ell
2^i}{d}\right)\right\}.
$$
The factor of $d\log L$ is an upper bound on the number of choices
for $i$; the factor of $2^{\ell i}$ is for the choice of the
$c_j$'s; and the factor ${L^d \choose \ell}$ is for the choice of
the $\ell$ (distinct) $v_j$'s. By (\ref{inner.property.1}) and the
second inequality in (\ref{inner.property.2}) we have (for $d$
sufficiently large)
$$
\eqalign{
2^{\ell i}{L^d \choose \ell}  \leq & ~~  2^{\ell i}\left(\frac{L^d}{\ell}\right)^\ell \cr
 \leq & ~~ 2^{\ell i} \left(O\left(2^\frac{id}{d-1}i^2 d^{1/2}\right)\right)^\ell  \cr
 \leq & ~~  2^{4\ell i} \cr
 =  & ~~  \exp\left\{o\left(\frac{2^i}{d}\right)\right\},
}
$$
so that in fact the right-hand side of (\ref{largeeventriv}) is at
most
$$
d \log L ~\max_i
\exp\left\{-\Omega\left(\frac{2^i\ell}{d}\right)\right\}.
$$
Taking $\ell$ as small as possible we see that this is at most
$$
d \log L ~\max_i \exp\left\{-\Omega\left(\frac{2^iL^d}{d
2^{\frac{id}{d-1}}i^2d^{1/2}}\right)\right\}
$$
and taking $i$ as large as possible we see that this is at most
$\exp\{-4L^{d-1}/d^4 \log^2 L\}$. Putting these observation together
we obtain (\ref{large3}).

\section{Proof of Lemma \ref{lem-volume.bounds.from.gk}}
\label{sec-proofs}

Much of what follows is modified from \cite{Galvin} and \cite{GKRS}.
Our strategy is as follows. Let $\ul{p}=(c_1,v_1, \ldots, c_\ell,
v_\ell)$ be given. Set $\ul{p'}=(c_2 ,v_2, \ldots, c_\ell, v_\ell)$.
We will show
\begin{equation} \label{induction}
\frac{\pi_3(\cC_3(\ul{p}))}{\pi_3(\cC_3(\ul{p'}))} \leq
\exp\left\{-\Omega\left(\frac{c_1}{d}\right)\right\}
\end{equation}
from which the lemma follows by a telescoping product. To obtain
(\ref{induction}) we define a one-to-many map $\varphi$ from
$\cC_3(\ul{p})$ to $\cC_3(\ul{p'})$. We then define a flow
$\nu:\cC_3(\ul{p}) \times \cC_3(\ul{p'}) \rightarrow [0,\infty)$
supported on pairs $(\chi,\chi')$ with $\chi' \in \varphi(\chi)$
satisfying
\begin{equation}\label{eq-flow.out}
\forall \chi \in \cC_3(\ul{p}),  \sum_{\chi' \in \varphi(\chi)}
\nu(\chi,\chi') =1
\end{equation}
and
\begin{equation}\label{eq-flow.in}
\forall \chi' \in \cC_3(\ul{p'}), \sum_{\chi \in
\varphi^{-1}(\chi')} \nu(\chi,\chi') \leq
\exp\left\{-\Omega\left(\frac{c_1}{d}\right)\right\}.
\end{equation}
This easily gives (\ref{induction}).

For each $s \in \{\pm 1, \ldots, \pm d\}$, define $\sigma_s$, the
{\em shift in direction $s$}, by $\sigma_s(x)=x+e_s$, where $e_s$ is
the $s$th standard basis vector if $s>0$ and $e_s=-e_{-s}$ if $s<0$.
For $X \subseteq V$ write $\sigma_s(X)$ for $\{\sigma_s(x):x \in
X\}$. For $\grg \in \cW$ set $W^s = \{x \in
\partial_{int}W:\sigma_{-s}(x) \not \in W\}$.

Let $\chi \in \cC_3(\ul{p})$ be given. Arbitrarily pick $\grg \in
\Gamma(I) \cap \cW(c_1,v_1)$ and set $W=\ant \grg$. Write $f$ for
the map from $\{0,1,2\}$ to $\{0,1,2\}$ that sends $0$ to $0$ and
transposes
$1$ and $2$. For each $s \in \{\pm 1, \ldots, \pm d\}$ and $S
\subseteq W$ define the function $\chi^s_S:V \rightarrow \{0,1,2\}$
by
$$
\chi^s_S(v) = \left\{ \begin{array}{ll}
                 0 & \mbox{if $v \in S$} \\
                 \chi(v) & \mbox{if $v \in (W^s\setminus S)
                 \cup (V \setminus W)$} \\
                 f(\chi(\sigma_{-s}(v))) & \mbox{if $v \in W \setminus
                 W^s$} \\
                      \end{array}
              \right.
$$
and set $\varphi_s(\chi)=\{\chi^s_S:S \subseteq W^s\}$.

\begin{claim} \label{claim-shift.works}
$\varphi_s(\chi) \subseteq \cC_3(\ul{p'})$.
\end{claim}

\noindent {\em Proof: }
An easy
case analysis verifies $\varphi_s(\chi) \subseteq \cC_3$.
We begin with the observation that the graph $\partial_{int} W \cup
\partial_{ext} W$ is bipartite with bipartition
$(\partial_{int} W, \partial_{ext} W)$. This follows from
(\ref{contour.prop.0}). By \eqref{contour.prop.1}, $I \cap
(\partial_{int} W \cup
\partial_{ext} W)=\emptyset$ and so for each component $U$ of $\partial_{int} W \cup
\partial_{ext} W$, $\chi$ is constant on $U \cap \partial_{int} W$
and on $U \cap \partial_{ext} W$ and in neither case does it take on
the value $0$.

Fix $S \subseteq W^s$. We show that if $\{u,v\}$ is an edge of
$T_{L,d}$ then $\chi^s_S(u) \neq \chi^s_S(v)$. We consider five
cases.

If $u, v \not \in W$ then $\chi^s_S(u)=\chi(u)$ and $\chi^s_S(v) =
\chi(v)$. But $\chi(u) \neq \chi(v)$, so $\chi^s_S(u) \neq
\chi^s_S(v)$ in this case.

If $u \in W$ and $v \not \in W$ then $\chi^s_S(v) = \chi(v)$ and
$\chi^s_S(u) \in \{0,\chi(u)\}$ (we will justify this in a moment).
Since $v \in \partial_{ext} W$ we have $\chi(v) \neq 0$ and we
cannot ever have $\chi(v)=\chi(u)$, so $\chi^s_S(u) \neq
\chi^s_S(v)$ in this case. To see that $\chi^s_S(u) \in
\{0,\chi(u)\}$, we consider subcases. If $u \in S$ then
$\chi^s_S(u)=0$. If $u \in W^s \setminus S$ then
$\chi^s_S(u)=\chi(u)$. Finally, if $u \in W \setminus W^s$ then
$\chi^s_S(u)=f(\chi(\gs_{-s}(u)))$; and $f(\chi(\gs_{-s}(u)))$ is
either $0$ or $\chi(u)$ depending on whether $\chi(\gs_{-s}(u))$
equals $0$ or $\chi(v)$ ($\chi(\gs_{-s}(u))$ cannot equal~$\chi(u)$).

If $u, v \in W\setminus W^s$ then $\chi^s_S(u)=f(\chi(\gs_{-s}(u)))$
and $\chi^s_S(v)=f(\chi(\gs_{-s}(v)))$. Since $f$ is a bijection and
$\chi(\gs_{-s}(u)) \neq \chi(\gs_{-s}(v))$ we have $\chi^s_S(u) \neq
\chi^s_S(v)$ in this case.

If $u \in W\setminus W^s$ and $v \in W^s \setminus S$ then
$\chi^s_S(u) \in \{0,\chi(u)\}$ (as in the second case above) and
$\chi^s_S(v)=\chi(v)$. Since $\chi(v) \neq 0$, we have $\chi^s_S(u)
\neq \chi^s_S(v)$.

Noting that it is not possible to have both $u, v \in W^s$, we
finally treat the case where $u \in W\setminus W^s$ and $v \in S$.
In this case $\chi^s_S(v)=\chi(v)=0$. Suppose (for a contradiction)
that $\chi^s_S(u)=0$. This can only happen if $\chi(\gs_{-s}(u))=0$.
If $\gs_{-s}(u)=v$, we have a contradiction immediately. Otherwise,
we have $\gs_{-s}(v) \not \in W$ and so (since
$\gs_{-s}(u)\gs_{-s}(v) \in E$) $\gs_{-s}(u) \in
\partial_{int} W$, also a contradiction.

This verifies $\varphi_s(\chi) \subseteq \cC_3$.
Because $\ant \grg$ is disjoint from the interiors of the remaining
cutsets in $\Gamma(I)$ and the operation that creates the elements
of $\varphi_s(\chi)$ only modifies $\chi$ inside $W$ it follows that
$\varphi_s(\chi) \subseteq \cC_3(\ul{p'})$. \qed

\begin{claim} \label{claim-unique.reconstruction}
Given $\chi' \in \varphi_s(\chi)$, $\chi$ can be uniquely
reconstructed from $W$ and $s$.
\end{claim}

\noindent {\em Proof: }Following \cite{GKRS}, we may reconstruct $\chi$ as follows.
$$
\chi(v) = \left\{ \begin{array}{ll}
                 \chi'(v) & \mbox{if $v \in V \setminus W$} \\
                 f(\chi'(\sigma_s(v))) & \mbox{if $v \in W$}. \\
                      \end{array}
              \right.
$$
\qed

We define the one-to-many map $\varphi$ from $\cC_3(\ul{p})$ to
$\cC_3(\ul{p'})$ by setting $\varphi(\chi)=\varphi_s(\chi)$ for a
particular direction $s$.
To define $\nu$ and $s$, we employ the notion of approximation also
used in \cite{GalvinKahn} and based on ideas introduced by
Sapozhenko in \cite{Sap}.
For $\grg \in \cW$, we say $A \subseteq V$ is an {\em approximation}
of $\grg$ if
$$
A^\cE \supseteq W^\cE ~~~\mbox{and}~~~A^\cO \subseteq W^\cO,
$$
$$
d_{A^\cO}(x) \geq 2d-\sqrt{d} ~~\mbox{for all $x \in A^\cE$}
$$
and
$$
d_{\cE \setminus A^\cE}(x) \geq 2d-\sqrt{d} ~~\mbox{for
all $y \in \cO \setminus A^\cO$},
$$
where $d_X(x)=|\partial x \cap X|$. Note that
from (\ref{contour.prop.0}) and (\ref{contour.prop.2}),
$W(\grg)$ is an approximation of $\grg$.

Before stating our main approximation lemma, which is a slight
modification of \cite[Lemma 2.18]{GalvinKahn}, it will be convenient
to further refine our partition of cutsets. To this end set
$$
\cW(w_e,w_o,v) = \left\{\grg:
\begin{array}{l}
\mbox{$\grg \in \Gamma(I)$ for some $\chi \in \cC_3^{even}$} \\
\mbox{with $|W^\cO|=w_o$, $|W^\cE|=w_e$}\\
\mbox{and $v \in W^\cE$}
\end{array}
\right\}.
$$
Note that
by (\ref{contour.prop.2}) we have
$|\grg| =2d(|W^\cO|-|W^\cE|)$ so $\cW(w_e,w_o,v) \subseteq
\cW((w_o-w_e)/2d, v)$.

\begin{lemma} \label{lem-k.comp.approx}
For each $w_e$, $w_o$ and $v$ there is a family $\cA(w_e,w_o,v)$
satisfying
$$
|\cA(w_e,w_o,v)| \leq
\exp\left\{O\left((w_o-w_e)d^{-\frac{1}{2}}\log^{\frac{3}{2}}d\right)\right\}
$$
and a map $\pi:\cW(w_e,w_o,v) \rightarrow \cA(w_e,w_o,v)$ such that
for each $\grg \in \cW(w_e,w_o,v)$, $\pi(\grg)$ is an approximation
for~$\grg$.
\end{lemma}

\noindent {\em Proof: } See \cite[Lemma 4.2]{Galvin}. \qed

We are now in a position to define $\nu$ and $s$. Our plan for each
fixed $\chi' \in \cC_3(\ul{p'})$ is to fix $w_e, w_o$ and $A \in
\cW(w_e, w_o,v)$
and to consider the contribution to the sum in (\ref{eq-flow.in})
from those $\chi \in \varphi^{-1}(\chi')$ with $\pi(\grg)=A$ (where
for each $\chi$, $\grg$ is a particular $\grg \in \Gamma(I) \cap
\cW(c_1,v_1)$).
We will try to define $\nu$ in such a way that each of these
individual contributions to (\ref{eq-flow.in}) is small; to succeed
in this endeavor we must first choose $s$ with care. To this end,
given $\grg \in \cW(w_e,w_o,v)$, set
$$
Q^\cE = A^\cE \cap \partial_{ext}(\cO \setminus A^\cO)
~~\mbox{and}~~Q^\cO = (\cO \setminus A^\cO) \cap
\partial_{ext} A^\cE,
$$
where $A = \pi(\grg)$ in the map guaranteed by Lemma
\ref{lem-k.comp.approx}. To motivate the introduction of $Q^\cE$ and
$Q^\cO$, note that for $\grg \in \pi^{-1}(A)$ we have
(by (\ref{contour.prop.0}) and (\ref{contour.prop.2}))
$$A^\cE\setminus Q^\cE \subseteq W^\cE,$$
$$\cE \setminus A^\cE \subseteq
\cE \setminus W^\cE,$$
$$A^\cO  \subseteq  W^\cO,$$
and
$$\cO \setminus (A^\cO \cup Q^\cO) \subseteq  \cO \setminus W^ \cO.$$
It follows that for each $\grg \in \pi^{-1}(A)$, $Q^\cE \cup Q^\cO$
contains all
vertices whose location in the partition
$T_{L,d} = W \cup \overline{W}$ is as yet unknown. We choose
$s(\chi)$ to be the smallest $s$ for which both of $|W^s| \geq
.8(w_o-w_e)$ and $|\sigma_s(Q^\cE) \cap Q^\cO| \leq 5|W^s|/\sqrt{d}$
hold.
This is the direction that minimizes the uncertainty
to be resolved when we attempt to reconstruct
$\chi$ from the partial information provided by $\chi' \in
\varphi^{-1}(\chi)$, $s$ and $A$.
(That such an $s$ exists is established in \cite[(49) and
(50)]{GalvinKahn} by an easy averaging argument).
Note that $s$ depends on $\grg$ but not~$I$.

Now for each $\chi \in \cC_3(\ul{p})$ let $\grg \in \Gamma(I)$ be a
particular cutset with $\grg \in \cW(c_1, v_1)$. Let $\varphi(\chi)$
be as defined before, with $s$ as specified above.
Define
$$C = W^s \cap A^\cO \cap \gs_s(Q^\cE)$$
and
$$D = W^s \setminus C,$$
and for each $\chi' \in \varphi(\chi)$ set
$$
\nu(\chi,\chi') = \left(\frac{1}{4}\right)^{|C \cap I(\chi')|}
\left(\frac{3}{4}\right)^{|C \setminus
I(\chi')|}\left(\frac{1}{2}\right)^{|D|}.
$$
Note that for $\chi \in \varphi^{-1}(\chi')$, $\nu(\chi,\chi')$
depends on $W$ but not on $\chi$ itself.

Since $C \cup D$ partitions $W$ we easily have (\ref{eq-flow.out}).
To obtain (\ref{induction}) and so (\ref{inq-bound.on.inner}) we
must establish (\ref{eq-flow.in}).

Fix $w_e$, $w_o$ such that $2d(w_o-w_e)=c_1$. Fix $A \in
\cA(w_e,w_o,v_1)$ and $s \in \{\pm 1, \ldots, \pm d\}$. For $\chi$
with $\grg \in \cW(w_e,w_o, v_1)$ write $\chi \sim_s A$ if it holds
that $\pi(\grg)=A$ and $s(\chi)=s$. We claim that with $A, s, w_o$
and $w_e$ fixed, for $\chi' \in \cC_3(\ul{p'})$
\begin{equation} \label{inq-black.box}
\sum \left\{\nu(\chi,\chi'):\chi \sim_s A,~\chi \in
\varphi^{-1}(\chi')\right\} \leq
\left(\frac{\sqrt{3}}{2}\right)^{w_o-w_e}.
\end{equation}
We could extract this directly from \cite{GKRS}, but for the
convenience of the reader we describe a proof below.

Write $\cC_3(\ul{p})(w_e,w_o,s,A,\chi')$ for the set of all $\chi
\in \cC_3(\ul{p})$ such that $W \in \cW(w_e,w_o,v_1)$,
$\pi(\grg)=A$, $s(\chi)=s$ and $\chi' \in \varphi(\chi)$ and set $U=
Q^\cE \cap \gs_{-s}(\chi')$. Say that a triple $(K,L,M)$ is {\em
good} for $\chi$ if it satisfies the following conditions.
$$
\mbox{$K \cup L \cup M$ is a minimal vertex cover of $Q^\cE \cup
Q^\cO$,}
$$
$$
\mbox{$K \subseteq Q^\cO$, $L \subseteq U$ and $M \subseteq Q^\cE
\setminus U$}
$$
and
$$
\mbox{$K=\partial_{ext}(U \setminus L)$}.
$$
We begin by establishing that $\chi \in
\cC_3(\ul{p})(w_e,w_o,s,A,\chi')$ always has a good triple.

\begin{lemma} \label{lem-core-exists}
For each $\chi \in \cC_3(\ul{p})(w_e,w_o,s,A,\chi')$ the triple
$$
(\hat{K}, \hat{L}, \hat{M}):=(W \cap Q^\cO, U \setminus W, (Q^\cE
\setminus U) \setminus W)
$$
is good for $\chi$.
\end{lemma}

\noindent {\em Proof: }\cite[around discussion of (54)]{GalvinKahn}.
\qed

In view of Lemma \ref{lem-core-exists} there is a triple $(K, L, M)$
that is good for $\chi$ and which has $|K|+|L|$ as small as
possible. Choose one such, say $(K_0(\chi), L_0(\chi), M_0(\chi))$.
Set $K'(\chi)=K_0\setminus \hat{K}$ and $L'(\chi)=L_0\setminus
\hat{L}$. Lemma \ref{lem-unlabeled_main_bound} below establishes an
upper bound on $\nu(\chi,\chi')$ in terms of $|K_0|$, $|L_0|$,
$|K'|$ and $|L'|$, and Lemma \ref{lem-reconstruction} shows that for
each choice of $K'$, $L'$ there is at most one $\chi$ contributing
to the sum in the lemma. These two lemmas combine to give
(\ref{inq-black.box}).
\begin{lemma} \label{lem-unlabeled_main_bound}
For each $\chi \in \cC_3(\ul{p})(w_e,w_o,s,A,\chi')$,
\begin{eqnarray*}
\nu(\chi,\chi') & \leq & \left(\frac{\sqrt{3}}{2}\right)^{w_o-w_e}
\frac{2^{|K_0|}}{3^{|K_0|+|L_0|}
2^{|K'|-|L'|}} \\
& := & B(K',L').
\end{eqnarray*}
\end{lemma}

\noindent {\em Proof: }We follow \cite[from just before (55) to just
after (60)]{GalvinKahn}, making superficial changes of notation.
\qed

The inequality in Lemma \ref{lem-unlabeled_main_bound} is the
$3$-coloring analogue of the main inequality of \cite{GalvinKahn}.
The key observation that makes this inequality useful is the
following.
\begin{lemma} \label{lem-reconstruction}
For each $w_e$, $w_o$, $s$, $A$, $\chi'$, $K'$ and $L'$, there is at
most one $\chi$ with $\chi \in \cC_3(\ul{p})(w_e,w_o,s,A,\chi')$,
$K'=K'(\chi)$ and $L'=L'(\chi)$.
\end{lemma}

\noindent {\em Proof: }In \cite[(56) and following]{GalvinKahn} it
is shown that $K'$ and $L'$ determine $W^\cO$ via
$$
\hat{K}= (K_0\setminus K') \cup (\partial_{ext} L' \cap Q^\cO)
$$
and so $W$ (via $W^\cE=\{v \in \cE:\partial v \subseteq W^\cO\}$).
But then by Claim \ref{claim-unique.reconstruction} $K'$ and $L'$
determine $\chi$. \qed

Lemmas \ref{lem-unlabeled_main_bound} and \ref{lem-reconstruction}
together now easily give (\ref{inq-black.box}):
\begin{eqnarray*}
\sum_{\chi \in \cC_3(\ul{p})(w_e,w_o,s,A,\chi')} \!\!\!\!\!
\nu(\chi,\chi') & \leq & \sum_{K' \subseteq K_0, ~L' \subseteq L_0}
\!\!\!B(K', L')
 \\
& \leq & \left(\frac{\sqrt{3}}{2}\right)^{w_o-w_e}.
\end{eqnarray*}

We have now almost reached (\ref{eq-flow.in}).
With the steps justified below we have that for each $\chi' \in
\cC_3(\ul{p'})$
\begin{eqnarray}
\sum_{\chi \in \varphi^{-1}(\chi')}  \nu (\chi,\chi')  & \leq &
\sum^{~~~~~\prime} \left\{\nu(\chi,\chi')\!: \!\begin{array}{l}\chi \sim_s A, \\
\chi \in
\varphi^{-1}(\chi')\end{array}\right\} \nonumber \\
& \leq & 2d c_1^{\frac{2d}{d-1}}|\cA(w_e, w_o,
v_1)|\left(\frac{\sqrt{3}}{2}\right)^{\frac{c_1}{2d}}
\label{overcount} \\
& \leq & 2d
c_1^{\frac{2d}{d-1}}\exp\left\{-\Omega\left(c_1/d\right)\right\}
\label{overcount2} \\
& \leq & \exp\left\{-\Omega\left(c_1/d\right)\right\},
\label{overcount3}
\end{eqnarray}
completing the proof of (\ref{eq-flow.in}). In the first inequality,
$\sum^\prime$ is over all choices of $w_e$, $w_o$, $s$ and  $A$. In
(\ref{overcount}), we note that there are $|\cA(w_e, w_o, v_1)|$
choices for
$A$, $2d$ choices for $s$ and
$c_1^{d/(d-1)}$ choices for each of $w_e$, $w_o$ (this is because
$c_1\geq (w_e+w_o)^{1-1/d}$, by (\ref{contour.prop.3})), and we
apply (\ref{inq-black.box})
to bound the summand. In (\ref{overcount2}) we use Lemma
\ref{lem-k.comp.approx}. Finally in (\ref{overcount3}) we use $c_1
\geq d^{1.9}$ (again by (\ref{contour.prop.3})) to bound
$2dc_1^{2d/(d-1)}=\exp\{o(c_1 /d)\}$.

\end{document}